\documentclass[preprint,12pt]{elsarticle}

\usepackage{amssymb,color,amsmath,textcomp,floatflt}

\newcommand{\eps}{\varepsilon}

\newtheorem{thm}{Theorem}
\newtheorem{cor}{Corollary}
\newtheorem{df}{Definition}
\newtheorem{lem}[thm]{Lemma}
\newdefinition{rem}{Remark}
\newdefinition{prop}{Proposition}
\newproof{pf}{Proof}

\begin{document}

\begin{frontmatter}



\title{A continuation principle for periodic BV-continuous state-dependent sweeping processes}







\author[mk]{Mikhail Kamenskii}
\author[om]{Oleg Makarenkov\corref{cor1}}
\ead{makarenkov@utdallas.edu}
\author[om]{Lakmi Niwanthi Wadippuli}

\cortext[cor1]{Corresponding author}

\address[mk]{Department of Mathematics, Voronezh State University, 394006 Voronezh, Russia}
\address[om]{Department of Mathematical Sciences, University of Texas at Dallas, 75080 Richardson, USA}

\sloppy

\begin{abstract}
We consider a Caratheodory differential equation with a state-dependent convex constraint that changes BV-continuously in time (a perturbed BV-continuous state-dependent sweeping processes). By setting up an appropriate catching-up algorithm we prove solvability of the initial value problem. Then, for sweeping processes with $T$-periodic right-hand-sides, we prove the existence of at least one $T$-periodic solution.  Finally, we further consider a $T$-periodic sweeping process which is close to an autonomous sweeping process with a constant constraint and prove the existence of a $T$-periodic solution specifically located near the boundary switched equilibrium of the autonomous sweeping process.
\end{abstract}






\end{frontmatter}



\section{Introduction} 

\noindent A variety of applications in elastoplasticity, economy, electrical circuits (see Adly et al \cite{adlynew} and references therein) lead to a  constraint differential equation 
\begin{equation}\label{sw11}
   -\dot x(t)\in N_{A(t)}(x(t))+f(t,x(t)),\qquad x\in E, 
\end{equation}
with a convex moving set $t\mapsto A(t)$ of just bounded variation (with respect to the Hausdorff metric). Here $E$ is a finite-dimensional vector space and $N_A(x)$ is a so-called normal cone defined for closed convex $A\subset E$ as 
\begin{equation}\label{NC}
  N_A(x)=\left\{\begin{array}{ll}\left\{\xi\in E:\left<\xi,c-x\right>\le 0,\ {\rm for\ any }\ c\in A\right\},& {\rm if}\ x\in A,\\
  \emptyset,& {\rm if}\ x\not\in A.
\end{array}\right.
\end{equation}
Whereas the case of Lipschitz $t\mapsto A(t)$ always leads (Edmond-Thibault \cite{edm2}) to the existence and uniqueness of an absolutely continuous solution $x(t)$ for any initial condition (under natural assumptions on $f$), the case where $t\mapsto A(t)$ is a convex-valued function of bounded variation doesn't ensure solvability of (\ref{sw11}) in the class of absolutely continuous functions. That is why an extended concept of the derivative (called Radon-Nikodym concept) is required in (\ref{sw11})  when the map $t\mapsto A(t)$ is a function of bounded variation, in which case equation (\ref{sw11}) is usually formulated in terms of differential measure 
$dx$ of BV continuous function $x$ and Lebesgue measure $dt$ as 
\begin{equation}\label{sw0}
   -dx\in N_{A(t)}(x)+f(t,x)dt,\quad x\in E.
\end{equation}
Existence and uniqueness of solutions to (\ref{sw0}) 
and 
 solvability of a periodic problem has been established by Castaing and Monteiro Marques in \cite{castaing}. The problem of existence and uniqueness of solutions in the unperturbed case ($f\equiv0$) was addressed in 
Moreau \cite{moreau}, Monteiro Marques \cite{26}, 
Valadier \cite{valadier}. Further state-independent extensions of (\ref{sw0}) were considered in  Adly et al \cite{adlynew}, Edmond-Thibault \cite{edm2},
Colombo and Monteiro Marques \cite{colombo}.

\vskip0.2cm

\noindent Motivated by applications to elastoplastic models with hardening and softening (see e.g. Gilles-Ulisse \cite{gilles}, Kunze \cite{kunze}), the present paper investigates the initial-value and periodic problems to the following state-dependent version of (\ref{sw0}): 
\begin{equation}\label{sw1}
   -dx\in N_{A+a(t)+c(x)}(x)+f(t,x)dt,\quad x\in E,
\end{equation}
where $a$ is a BV-continuous function and $c:E\mapsto E$ is a Lipschitz function.

\vskip0.2cm

\noindent To prove the existence of solutions to (\ref{sw1}) we introduce a new catching-up algorithm which allows to construct a sequence $\{x_n\}_{n\in\mathbb{N}}$ of approximations of the solution $x$ of (\ref{sw1}).  
 The existence of $T$-periodic solutions to (\ref{sw1}) is proved by establishing that 
\begin{equation}\label{d}
d(I-P^{n},Q)\not=0
\end{equation} for 
the Poincare maps $P^{n}$ of the $n$-th approximation of the catching-up scheme and suitable $Q\in E.$ 
Here $d(I-P^n,Q)$ is the topological degree of the map $P^n$ with respect to an open bounded set $Q$ (see Krasnoselskii-Zabreiko \cite{kz}).
After we get the existence of a fixed point for $P^{n}$ we pass to the limit as $n\to\infty$ on the respective $T$-periodic solutions of sweeping process (\ref{sw1}) and get the existence of a $T$-periodic solution to (\ref{sw1}) even though we don't know whether $\lim_{n\to\infty} P^{n}(x)$ exists or not.

\vskip0.2cm

\noindent The paper offers global and local sufficient conditions to ensure (\ref{d}). The global sufficient condition is based on  construction of such a convex set $Q$ which contains all possible values of the set $A+a(t)+c(x(t))$ for all possible solutions of (\ref{sw1}). In this way, we can show that $P^{n}(Q)\subset \overline Q$ for sufficiently large $n\in\mathbb{N},$ which ensures (\ref{d}). 

\vskip0.2cm

\noindent To design sufficient conditions that ensure the validity of (\ref{d}) in a desired region $Q$ (local sufficient conditions), we are no longer allowed to enlarge $Q$ as much as we want, so we have to seek for alternative deformations of (\ref{d}) that stick to the given region $Q$. We go here by a continuation approach and replace (\ref{sw1}) by a parameter dependent sweeping process 
\begin{equation}\label{sw1lam}
   -dx\in N_{A+a(t,\lambda)+c(x,\lambda)}(x)+f(t,x,\lambda)dt,\quad  x\in E, \ \lambda\in\mathbb{R}.
\end{equation}
 Accordingly, the relation (\ref{d}) gets replaced  by 
\begin{equation}\label{dlam}
d(I-P^{\lambda,n},Q)\not=0.
\end{equation} 
We, therefore, assume that (\ref{d}) corresponds to (\ref{dlam}) with some $\lambda=\lambda_1$ and prove the validity of (\ref{dlam}) for $\lambda=\lambda_1$ building upon some good properties of $P^{\lambda,n}$ for $\lambda=0$ combined with nondegenerate homotopy between $P^{\lambda_1,n}$ and $P^{0,n}.$ 

\vskip0.2cm

\noindent As for possible good properties of $P^{0,n}$ we offer both topological (Theorem~\ref{thm3}) and algebraic (Theorem~\ref{thm16})  conditions. The topological condition simply assumes that (\ref{dlam}) holds for $\lambda=0$, that leads to an analogue of standard continuation principles available for ordinary differential equations, see e.g. Capietto et al \cite{mawhin} and Kamenskii et al \cite{nachr}.

\vskip0.2cm

\noindent To obtain easily verifiable algebraic conditions that ensure the validity of (\ref{dlam}) for $\lambda=0$, this paper takes a straightforward route and simply offers sufficient conditions for asymptotic stability of a point $x_0$ of the target set $Q.$ Such an approach is based on the fact that the topological degree of a Poincare map in the neighborhood of an asymptotically stable fixed point equals 1.  However, just assuming that $x_0$ is an asymptotically stable equilibrium of (\ref{sw1lam}) with $\lambda=0$ is not of interest because it leads to periodic solutions that  don't interact with the boundary of the constraint of (\ref{sw1lam}) when $\lambda>0.$ Such periodic solutions will simply be solutions of the differential equation 
\begin{equation}\label{ode}
   -dx=f(t,x,\lambda)dt,\quad  x\in E, \ \lambda\in\mathbb{R}.
\end{equation}
That is why a non-equilibrium concept of an asymptotically stable point $x_0$ is required  to design periodic solutions of (\ref{sw1lam}) which are intrinsically sweeping (i.e. interact with the boundary of the constraint of (\ref{sw1lam})).

\vskip0.2cm

\noindent The required concept of asymptotically stable point $x_0$ has been recently developed in Kamenskii-Makarenkov \cite{kamenskii} based on the notion of switched boundary equilibrium well known  in control theory (see e.g. Bolzern-Spinelli \cite{bol}). To introduce the concept of switched boundary equilibrium for sweeping process (\ref{sw1lam}) at $\lambda=0$, we will assume that, at $\lambda=0$, sweeping process (\ref{sw1lam}) takes the form 
\begin{equation}\label{aaa}-\dot x\in N_{A}(x)+f_0(x), \qquad x\in E,
\end{equation}
where $A$ is just a constant convex closed bounded set and $f_0$ is Lipschitz continuous. Following  Kamenskii-Makarenkov \cite{kamenskii}, $x_0$ is a switched boundary equilibrium of (\ref{aaa}), if $f_0(x_0)$ is normal to the boundary $\partial A$ and if $f_0(x_0)$ points inwards $A$. To prove the validity of (\ref{dlam}) for $\lambda=0$ this paper, therefore, offers sufficient conditions for asymptotic stability of switched boundary equilibrium $x_0$ of (\ref{aaa}).



\vskip0.2cm

\noindent The paper is organized as follows. In the next section of the paper we introduce a formal definition of sweeping process (\ref{sw1lam}) following  Castaing and Monteiro Marques \cite{castaing}. The fundamental result of the paper (Theorem~\ref{thm1}) on solvability of the initial-value problem for (\ref{sw1}) (and, therefore, for (\ref{sw1lam}) too) is formulated in Section~\ref{3}. In the same section we introduce our concept of {\it generalized initial condition} that we repeatedly use in the paper later and which allows us to consider solutions of (\ref{sw1}) with initial conditions outside $A+a(t)+c(x)$. We simply say that $x(t)$ is a solution of (\ref{sw1}) with initial condition $q\in E$, if $x(0)$ is the solution of $x(0)={\rm proj}(q,A+a(0)+c(x(0)))$ (which has a unique solution $x(0)=V^0(q)$ according to Lemma~\ref{lem7}).

\vskip0.2cm

\noindent Sections~\ref{secperiodic} and \ref{4} contain formulations of our results on the existence of $T$-periodic solutions to (\ref{sw1}). 
 Section~\ref{secperiodic} offers a theorem (Theorem~\ref{thmexist}) saying that any $T$-periodic state-dependent sweeping process (\ref{sw1}) always admits at least one $T$-periodic solution, if the right-hand-sides of (\ref{sw1}) are $T$-periodic. Remarkably, the theorem doesn't assume uniqueness or continuous dependence of solutions of (\ref{sw1}) on initial conditions.
 

\vskip0.2cm

\noindent Abstract results on { continuation} of $T$-periodic solutions to (\ref{sw1lam}) are presented in Section~\ref{4}. We assume that for $\lambda=0$ the sweeping process (\ref{sw1lam}) admits a Poincare map $P^0$ (over time $T$) and formulate (Theorem~\ref{thm3}) a standard continuation principle: {\it if the topological degree $d(I-P^0\circ V^0,Q)\not=0$ for some open bounded set $Q\subset E$ and if none points of the boundary of $Q$ are initial conditions of $T$-periodic solutions of sweeping process (\ref{sw1lam}) for any $\lambda\in[0,\lambda_1]$ (non-degenerate deformation), then, for any $\lambda\in[0,\lambda_1]$, sweeping process (\ref{sw1lam}) admits a $T$-periodic solution $x$.} A result on the existence of $\lambda_1>0$ such that the non-degenerate deformation assumption of Theorem~\ref{thm3} holds is also presented (Theorem~\ref{thm4}) in Section~\ref{4}.  

\vskip0.2cm

\noindent Section~\ref{5} contains  proofs of  Theorems~\ref{thm1}-\ref{thm4}. The proof of the existence of solutions is based on introducing (section~\ref{sec62}) a new implicit catching-up scheme (\ref{c1})-(\ref{c4}), which in turn relies on the following two ideas: {\bf (i)}  Castaing and Monteiro Marques   change of the variables \cite[Theorem~4.1]{castaing} that  converts (section~\ref{equivalent}) the perturbed sweeping process (\ref{sw1lam}) with differential measure $dx$ into a non-perturbed sweeping process 
(\ref{sw3}) for the derivative 
    $\frac{du}{|du|}$ with respect to the variation measure $|du|$ of $du$; {\bf (ii)} Kunze and Monteiro Marques lemma (\cite[Lemma~7]{kunze}) to resolve (Lemma~\ref{lem7}) the implicit catching-up scheme (\ref{c1})-(\ref{c4}) with respect to the implicit variable. 
Furthermore, our Lemma~\ref{lem7} extends \cite[Lemma~7]{kunze} by proving continuous dependence of scheme (\ref{c1})-(\ref{c4}) on initial condition, that gave us continuity of Poincare maps $P^{\lambda,n}$ (section~\ref{poincaremapsec}).    
    The convergence of the scheme (\ref{c1})-(\ref{c4}) is established in section~\ref{sec64} where we prove (Lemma~\ref{lemconverge})  convergence of the approximations $\{u_n\}_{n\in\mathbb{N}}$ of solution $u$ of (\ref{sw3}) and then prove (Lemma~\ref{lemkam1}) convergence of the respective approximations $\{x_n\}_{n\in\mathbb{N}}$  of solution $x$ of sweeping process (\ref{sw1lam}). In other words, Lemma~\ref{lemkam1} states that 
the change of the variables of     
    Castaing and Monteiro Marques  \cite[Theorem~4.1]{castaing} is continuous with respect to time-discretization.
    Finally, a result by     
    Monteiro Marques \cite[p.~15-16]{26} (which is also Proposition~6 in 
Valadier \cite{valadier}) is used to prove (Theorem~\ref{thmexistence} of section~\ref{sec65}) 
    that the limit of catching-up scheme (\ref{c1})-(\ref{c4}) is a solution of (\ref{sw1lam}). 

\vskip0.2cm

\noindent Section~\ref{7} is devoted to establishing conditions for the validity of (\ref{dlam}) at $\lambda=0$ in a neighborhood $Q$ of a switched boundary equilibrium $x_0$. Specifically, as mentioned earlier, we assume that, for $\lambda=0$ sweeping process (\ref{sw1lam}) takes the form (\ref{aaa}) and discover conditions for asymptotic stability of $x_0\in\partial A.$ 
In particular, in section~\ref{7} we extend the two-dimensional approach of Makarenkov and Niwanthi Wadippuli \cite{ijbc} and derive a differential equation of sliding motion along $\partial A$, for which $x_0$ is a regular equilibrium whose stability can be investigated (Theorem~\ref{stabslide}) over the eigenvalues of the respective linearization. Assuming that the real parts of these eigenvalues are negative we conclude that $d(I-P^0\circ V_0,Q)=1$ and establish (Theorem~\ref{thm16}) the existence of $T$-periodic solutions near $x_0$ for all BV-continuous state-dependent sweeping processes (\ref{sw1lam}) that approaches  (\ref{aaa}) when $\lambda\to 0.$

\vskip0.2cm

\noindent Conclusions and Acknowledgments sections conclude the paper.

\section{Definition of solution}

\noindent In what follows, $\mathcal{B}([0,T])$ is the family of Borel subsets of $[0,T].$ A {\it Borel vector measure} on $[0,T]$ is a map $\mu:\mathcal{B}([0,T])\to E$ such that $\mu\left(\cup_{n=1}^\infty B_n\right)=\Sigma_{n=1}^\infty\mu(B_n)$ for any sequence $\{B_n\}_{n=1}^\infty$ of mutually disjoint elements of $\mathcal{B}([0,T]),$ see Recupero \cite[\S2.4]{rec} or Dinculeanu \cite[Definition~1, \S\hskip-0.15cm~III.14.4, p.~297]{Dinculeanu}. 

\vskip0.2cm

\noindent According to Dinculeanu \cite[Theorem~1, \S\hskip-0.15cm~III.17.2, p.~358]{Dinculeanu} (see also  Recupero \cite{rec}), any BV-continuous function $x:[0,T]\to E$ admits a unique vector measure of bounded variation $dx:\mathcal{B}([0,T])\to E$ (called {\it Stieltjes measure} in \cite{Dinculeanu}) such that for every $0<t_1<t_2<T$ we have
$$
  \begin{array}{ll}
     dx((t_1,t_2))=x(t_2^-)-x(t_1^+), &  dx([t_1,t_2])=x(t_2^+)-x(t_1^-),\\
     dx([t_1,t_2))=x(t_2^-)-x(t_1^-), & dx((t_1,t_2])=x(t_2^+)-x(t_1^+),
  \end{array}
$$
where 
$$
  x(t^-)=\lim_{\tau\to t^-}x(\tau),\quad
  x(t^+)=\lim_{\tau\to t^+}x(\tau),\quad  0<t<T.	
$$

\noindent A vector Borel measure $d\mu$ is called continuous with respect to a scalar Borel measure  $d\nu$ (or simply $d\nu$-continuous), if $\lim_{\nu(D)\to 0}\mu(D)=0$, see Diestel-Uhl \cite[p.~11]{diestel}. If a vector measure $d\mu$ is $d\nu$-continuous then, according to Radon-Nikodym Theorem \cite[p.~59]{diestel} there is a $d\nu$-integrable function $g:[0,T]\mapsto E$ such that
$$
   d\mu(D)=\int_D g\hskip0.05cm d\nu,\quad\mbox{for all}\ D\in\mathcal{B}([0,T]). 
$$
In this case, the function $g$ is called Radon-Nikodym derivative of $d\mu$ with respect to $d\nu$ (or density) and is denoted by $\dfrac{d\mu}{d\nu}.$ Furthermore, according to Moreau-Valadier \cite[Proposition~1]{mv} (see also Valadier \cite[Theorem~3]{valadier}), the Radon-Nikodym derivative $\dfrac{d\mu}{d\nu}$ can be computed as
$$
  \dfrac{d\mu}{d\nu}(t)=\lim\limits_{\eps\to 0,\hskip0.05cm \eps>0}\dfrac{d\mu([t,t+\eps])}{d\nu([t,t+\eps])}, \quad d\nu-a.e.\ on\ [0,T].
$$

\vskip0.2cm

\noindent We will use the following definition of the solution of (\ref{sw1}) (Castaing and Monteiro Marques \cite[\S1]{castaing}).
\begin{df}\label{df2}  A BV continuous function $x$ is called a solution of (\ref{sw1}), if there exists a finite measure $d\nu$ for which both  differential measure $dx$ and Lebesgue measure $dt$ are $d\nu$-continuous, and such that
$$
  -\frac{dx}{d\nu}(t)\in N_{A+a(t)+c(x(t))}(x(t))+f(t,x(t))\frac{dt}{d\nu}(t),\quad d\nu-a.e.\ on\ [0,T].
$$
\end{df}

\section{Existence of solutions}\label{3}


\noindent It is customary (see \cite[Theorem~6]{kunze2}) to assume that the initial condition $q$ of sweeping process (\ref{sw1}) satisfies
\begin{equation}\label{ic}
   q\in A+a(0)+c(q).
\end{equation}
However, it will be convenient for our analysis to define solutions of (\ref{sw1}) for any initial condition $q\in E,$ that we will term a {\it generalized initial condition}. We take advantage of the fact, that for contracting map $c$, the equation
$$
  v={\rm proj}(q,A+a(0)+c(v))
$$
always has a solution $v=V(q)$ (see Lemma~\ref{lem7}) and $V\in C^0(E,E).$ 

\vskip0.2cm

\noindent The core of this paper is the following Theorem~\ref{thm1} on the existence of solutions to (\ref{sw1}). As itself, the theorem won't loss anything by dropping the generalized initial condition concept. However, considering the generalized initial conditions will be convenient for applications of Theorem~\ref{thm1} to the problem of the occurrence of periodic solutions from a boundary equilibrium, that we consider in this paper later (Theorem~\ref{thm16}).

\begin{thm} \label{thm1} Assume that $A\subset E$ is a nonempty closed convex bounded set, $a:[0,T]\to E$ is BV-continuous on $[0,T]$,  $x\mapsto c(x)$ is globally Lipschitz with Lipschitz constant $0<L_2<1$, and $(t,x)\mapsto f(t,x)$ is Caratheodory in $(t,x)$  with respect to Lebesgue measure and globally Lipschitz in $x.$ Then, for any generalized initial condition $q\in E$,
the sweeping process (\ref{sw1}) admits a solution, defined on $[0,T],$ with the initial condition $x(0)=V(q)$. In particular, sweeping process (\ref{sw1}) admits a solution on $[0,T],$ for any initial condition $x(0)=q$, where $q$ satisfies  (\ref{ic}).
\end{thm}



\section{Global existence of periodic solutions}\label{secperiodic}

\noindent In this section we offer a result saying that, under the conditions of Theorem~\ref{thm1}, sweeping process (\ref{sw1}) always has a periodic solution, if the right-hand-sides are $T$-periodic.

\begin{thm}\label{thmexist} Assume that conditions of theorem~\ref{thm1} hold and let $L_2\in(0,1)$ be the Lipschitz constant of $c$ as introduced in theorem~\ref{thm1}. Denoting by $\xi\in E$  the unique solution of $c(\xi)=\xi,$ consider the set
$$\Omega=\bigcup\limits_{t\in[0,T]}\Omega_t,\quad \Omega_t=\bigcup\limits_{b\in A(t)}\left\{x:\|x-\xi\|<\dfrac{\|b\|}{1-L}\right\}.$$
 Then sweeping process (\ref{sw1}) admits a solution $t\mapsto x(t)$ such that
\begin{equation}\label{Tper}
   x(T)=x(0)\in \overline\Omega.
\end{equation}
In particular, $t\mapsto x(t)$ is a $T$-periodic solution of (\ref{sw1}), if both $t\mapsto a(t)$ and $t\mapsto f(t,x)$ are $T$-periodic.
\end{thm}

\begin{rem} Throughout the paper we prefer to work with functions defined on $[0,T]$ only. When saying $t\mapsto x(t)$ is a $T$-periodic solution of (\ref{sw1}), we mean  that $t\mapsto x(t)$ becomes a $T$-periodic solution after all functions are extended to $\mathbb{R}$ by $T$-periodicity.
\end{rem}

\section{Continuation of periodic solutions}\label{4}

\noindent This section considers a  $\lambda$-dependent sweeping process (\ref{sw1lam}) for measures $dx$ and $dt$, and discovers how the existence of periodic solutions for $\lambda>0$ can be concluded from an appropriate knowledge about (\ref{sw1lam}) at $\lambda=0.$

\vskip0.2cm

\noindent We will assume that BV-continuity of $a$ of Theorem~\ref{thm1} holds uniformly with respect to $\lambda,$ i.e.  
\begin{equation}\label{A1_new}
\begin{array}{l}
   {\rm var}(a(\cdot,\lambda),[s,t])\le {\rm var}(\bar a,[s,t]),\quad \lambda\in[0,1],\\
\mbox {where }\bar a:[0,T]\to\mathbb{R}\mbox{ is a BV continuous function}.
\end{array}
\end{equation}


\noindent The map $V^\lambda$ for (\ref{sw1lam}) now depends on the parameter $\lambda$ and is defined as the unique solution (according to Lemma~\ref{lem7}) of the equation
$$
  v={\rm proj}(q,A+a(0,\lambda)+c(v,\lambda)).
$$

\noindent We will call sweeping process (\ref{sw1lam}) $T$-periodic, if 
$$
   a(t+T,\lambda)\equiv a(t,\lambda),\quad f(t+T,x,\lambda)\equiv f(t,x,\lambda).
$$

\noindent In what follows, $d(I-\bar P,Q)$ is the topological degree of the vector field $I-\bar P$ on an open bounded set $Q\subset E,$ see e.g. Krasnoselskii-Zabreiko \cite{kz}.

\begin{thm} \label{thm3} Assume that $T$-periodic sweeping process (\ref{sw1lam}) possesses the following regularity:
\begin{itemize}
\item[{\rm I)}] The set $A\subset E$ is nonempty, convex, closed, and bounded. The function $a$ satisfies (\ref{A1_new}). The function $x\mapsto c(x,\lambda)$ is globally Lipschitz with Lipschitz constant $0<L_2<1.$ The function $(t,x)\mapsto f(t,x,\lambda)$ is Caratheodory in $(t,x)$ with respect to Lebesgue measure and globally Lipschitz in $x,$ and both the Lipschitz constants are independent of $\lambda\in[0,1].$ Furthermore, $a$, $c$, and $f$ are continuous in $\lambda\in[0,1]$ uniformly with respect to $t\in[0,T]$ and $x\in E.$ 
\end{itemize}
Assume, that the existence of a $T$-periodic solution for $\lambda=0$ is given in the following extended way: 
\begin{itemize}
\item[{\rm II)}] There exists an open bounded $Q\subset E$ such that, when $\lambda=0$, the solution of (\ref{sw1lam}) is unique for any initial condition $x(0)\in V^0(\overline Q)$, none of the elements of $V^0(\partial Q)$ are initial conditions of $T$-periodic solutions of (\ref{sw1lam}) with $\lambda=0,$ and for the Poincare map $P^0$ of (\ref{sw1lam}) with $\lambda=0$ one has
$$
  d(I-P^0\circ V^0,Q)\not=0.
$$
\end{itemize}
Finally, assume the following homotopy  through $\lambda\in[0,\lambda_1]$: 
\begin{itemize}
\item[{\rm III)}] There exists $\lambda_1 \in (0,1]$ such that sweeping process
(\ref{sw1lam}) doesn't have periodic solutions $x$ with initial condition $x(0)\in V^\lambda(\partial Q),$ $\lambda\in[0,\lambda_1].$
\end{itemize}
 Then, for any $\lambda\in[0,\lambda_1]$, sweeping process (\ref{sw1lam}) admits a $T$-periodic solution $x$ with the initial condition $x(0)\in V^\lambda(Q).$
\end{thm}

\noindent Note, for $\lambda>0$, we don't know whether or not the solutions of sweeping process (\ref{sw1lam}) are uniquely defined by the initial condition or depend continuously on $\lambda$. That is why the statement of the theorem is not a direct consequence of II) as it usually happens  in topological degree based existence results. In particular, we cannot establish any type of continuity of solutions as $\lambda\to 0$. That is why the next theorem is not a direct consequence of Theorem~\ref{thm3}.

\begin{thm}\label{thm4}
Assume that sweeping process (\ref{sw1lam}) is $T$-periodic.
Assume that  conditions I) and II) of Theorem~\ref{thm3} hold. Then, there exists $\lambda_1>0$ such that  condition~III) of Theorem~\ref{thm3} holds, and, therefore, for any $\lambda\in[0,\lambda_1]$, sweeping process (\ref{sw1lam}) admits a $T$-periodic solution $x$ with the initial condition $x(0)\in V^\lambda(Q).$
\end{thm}

\section{The catching-up algorithm and proofs of the abstract existence results} \label{5}

\subsection{An equivalent non-perturbed formulation of the initial perturbed sweeping process}\label{equivalent}

\noindent Recall, that for a BV-continuous function $u:[0,T]\to E,$ the {\it variation measure} $|du|$ (also called {\it modulus measure}) is defined, for any $D\in \mathcal{B}([0,T]),$  as (see Diestel-Uhl \cite[Definition~4, p.~2]{diestel}, Recupero \cite[\S2.4]{rec})
\begin{eqnarray*}
&&\hskip-0.5cm |du|(D)=\\
&&\hskip-0.5cm =  \sup\left\{\sum_{n=1}^\infty\|u(D_n)\|:D=\bigcup_{n=1}^\infty D_n,\ D_n\in\mathcal{B}([0,T]),\ D_i\cap D_j=\emptyset\ {\rm if}\ i\not=j\right\}. 
\end{eqnarray*}
For a BV-continuous function $u:[0,T]\to\mathbb{R}$, the differential measure $du$
is always $|du|$-continuous (it follows e.g. from Diestel-Uhl \cite[Theorem~1, p.~10]{diestel}), i.e. a $|du|$-integrable  density 
$\dfrac{du}{|du|}$ is well defined. Moreover, according to Castaing and Monteiro Marques \cite[Theorem~4.1]{castaing}, if $x$ is a solution of the perturbed sweeping process (\ref{sw1lam}), then the $BV$ continuous function $u$ defined by  
\begin{equation}\label{xu}
   u(t)=x(t)+\int_0^t f(\tau,x(\tau))d\tau
\end{equation}
is a solution to the non-perturbed sweeping process
\begin{equation}\label{sw3}
    -\frac{du}{|du|}(t)\in N_{A+a(t,\lambda)+c(x(t),\lambda)+\int_0^t f(\tau,x(\tau),\lambda)d\tau}(u(t)),\quad |du|-a.e.\ on\ [0,T].
\end{equation}

\begin{lem}\label{lemkam} Assume that $(t,x,\lambda)\mapsto f (t,x,\lambda)$  is Caratheodory in $(t,x)$ with respect to Lebesgue measure and is globally Lipschitz in $x$ with Lipschitz constant independent of $t\in[0,T]$ and $\lambda\in[0,1].$  Then, for any BV continuous $u:[0,T]\to E$, the equation (\ref{xu})
admits a unique BV continuous solution $x:[0,T]\to E$.
\end{lem}

\noindent Lemma~\ref{lemkam} is a direct consequence of Lemma~\ref{lemkam1} that we prove below.

\vskip0.2cm

\noindent Combining  \cite[Theorem~4.1]{castaing} and Lemma~\ref{lemkam}, we can formulate the following equivalent definition of the solution of (\ref{sw1lam}).

\begin{df}\label{df1}  A BV continuous function $x$ is called a solution of perturbed sweeping process (\ref{sw1lam}), if the function $u$ given by (\ref{xu}) is a solution of the non-perturbed sweeping process  (\ref{sw3}). 
\end{df}

\subsection{The catching-up algorithm}\label{sec62}

\noindent For each fixed $n\in\mathbb{N},$ we partition $[0,T]$ into smaller intervals by the points $\{t_0,t_1,...,t_n\}\subset[0,T]$ defined by  $$\ t_0=0,\ t_n=T,\ t_{i+1}-t_i=\frac{T}{n}, \quad i\in\overline{1,n}.$$ In what follows, we fix some initial condition 
$$
   x(0)=u(0)=q,
$$
where $q$ satisfies 
\begin{equation}\label{q0lam}
     q\in A+a(0,\lambda)+c(q,\lambda),
\end{equation}
and use the ideas of Definition~\ref{df1} in order to construct pieceiwise-linear  functions $u_n$ and $x_n$ (linear on each $[t_{i},t_{i+1}]$) that  serve as approximations of the solutions $u$ and $x$ of Definition~\ref{df1}. The construction will be implemented iteratively through the intervals $[t_i,t_{i+1}]$ starting from $i=0,$ and moving towards $i=n-1.$ 

\vskip0.2cm

\noindent Denoting
$$
   u_n(0)=q,\ x_n(0)=q,\ u_i^n=u_n(t_i),\ x_i^n=x_n(t_i),\quad i\in\overline{0,n},
$$
we apply the implicit iterative scheme
\begin{eqnarray}
   u_{i+1}^n&=&{\rm proj}  \left[u_i^n,A+a(t_{i+1},\lambda)+c\left(u_{i+1}^n-\int\limits_0^{t_i}f(\tau,x_n(\tau),\lambda)d\tau, \lambda \right)\right.\nonumber\\
&& \left.+\int\limits_0^{t_i}f(\tau,x_n(\tau),\lambda)d\tau\right],\label{c1}\\
x^n_{i+1}&=&u^n_{i+1}-\int\limits_0^{t_i}f(\tau,x_n(\tau),\lambda)d\tau,\label{c2}\\
  	 u_n(t)&=& u^n_i + \dfrac{t-t_i}{t_{i+1}-t_i} (u^n_{i+1}-u^n_i),\quad t\in [t_i,t_{i+1}],\label{c3}\\
x_n(t)&=& x^n_i + \dfrac{t-t_i}{t_{i+1}-t_i} (x^n_{i+1}-x^n_i),\quad  t\in [t_i,t_{i+1}],\label{c4}
\end{eqnarray}
successively from $i=0$ to $i=n-1$. Next lemma uses the idea of the implicit scheme of Kunze and Monteiro Marques (\cite[Lemma~7]{kunze2})
and it
proves that  for each $i\in\overline{0,n-1}$ we can extend the definition of $u_n$ and $x_n$ from $[0,t_i]$ to $[0,t_{i+1}]$ according to (\ref{c1})-(\ref{c4}).

\begin{lem}\label{lem7} Consider a set-valued function 
$$
  C(s_1,s_2,u,\xi)=A+\tilde a(s_1, \xi)+\tilde c(s_2,u,\xi),\quad s_1,s_2\in[0,T],\ u\in E,\ \xi\in W,
$$
where $A\subset E$ is a nonempty closed convex bounded set,  $\tilde a:\mathbb{R}\times W\to E,$ $\tilde c:\mathbb{R}\times E\times W\to E$, and $W$ is a finite dimensional Euclidean space. Assume that 
$$
 \begin{array}{l} 
    {\rm var}(\tilde a(\cdot,\xi),[s,t]) \le{\rm var}(\bar a,[s,t]),\quad \xi\in W,\\
  \mbox {where }\bar a:[0,T]\to\mathbb{R}\mbox{ is a BV continuous function},
 \end{array}
$$
and  $(s,\xi)\to \tilde a(s,\xi)$ is continuous in $\xi\in W$ uniformly in $s\in[0,T].$ Assume that $(s,u,\xi)\mapsto \tilde c(s,u,\xi)$ is continuous in $\xi\in W$ uniformly in $(s,u)\in[0,T]\times E$ and satisfies the Lipschitz condition
\begin{eqnarray*} 
&&\lVert \tilde c(s,u,\xi) - \tilde c(t,v,\xi) \rVert \leq L_1 |s-t| +L_2 \lVert u-v \rVert,\\ && \qquad\qquad\mbox{for any }s,t\in[0,T],\ u,v\in E,\ \xi\in W, 
\end{eqnarray*}
 with $L_1>0$ and $L_2\in(0,1).$
	Then, for any $\tau_1,\tau_2,s_1,s_2\in[0,T]$ and any $u\in E$
 there exists an unique $v=v(\tau_1,\tau_2,s_1,s_2,u,\xi)$ such that 
\begin{equation}\label{reqprop}
v \in C(\tau_1,\tau_2,v,\xi)\quad{\rm and}\quad
v={\rm proj} (u,C(\tau_1,\tau_2,v,\xi)).
\end{equation} Moreover, $v\in C^0([0,T]\times[0,T]\times[0,T]\times[0,T]\times E\times W,E)$. If, in addition,  $$
u\in C(s_1,s_2,u,\xi),
$$ then
\begin{equation}\label{estimate}
\lVert v-u \rVert \leq\dfrac{{\rm var}(\bar a,{[s_1,\tau_1]}) +L_1 |\tau_2-s_2|}{1-L_2}.
\end{equation}
\end{lem}

\noindent The following key estimate is required for the proof of Lemma~\ref{lem7}.

\begin{lem}\label{extra1}
  Let  $C$ be a convex set of $E$.  Then, for any vectors $u,c\in E$, 
  $$ \|{\rm proj}(u,C)-{\rm proj}(u,C+c)\| \le \|c\|.$$ 
\end{lem}

\noindent {\bf Proof.} From the definition of projections $v_1={\rm proj}(u,C)$ and  $v_2={\rm proj}(u,C+c)$ we have (see e.g. Kunze and Monteiro Marques \cite[\S2]{kunze2})
\begin{equation}\label{concl}
u-v_1\in N_{C}(v_1)\quad {\rm and}\quad u-v_2\in N_{C+c}(v_2).
\end{equation}
Since 
$v_2-c\in C$ and $v_1+c\in C+c$, we conclude from (\ref{concl}) that 
$$
\left< u-v_1, v_2-c-v_1\right>\leqslant 0\quad{\rm and}\quad \left< u-v_2, v_1+c-v_2\right>\leqslant 0,
$$
or, rearranging the terms,
$$
\left< v_1-u, v_1-v_2\right>\leqslant \left< u-v_1,c\right> \quad {\rm and}\quad \left< u-v_2,v_1-v_2\right>\leqslant \left< v_2-u,c\right>.
$$
Finally, we add both inequalities together and get 
\[
\left< v_1-v_2,v_1-v_2\right>\leqslant \left< v_2-v_1,c\right> \leqslant\|v_1-v_2\|\cdot\|c\|, 
\]
which implies the statement. \qed

\vskip0.2cm

\noindent {\bf Proof of Lemma~\ref{lem7}.}	{\bf Step 1.} {\it The existence of $v(\tau_1,\tau_2,s_1,s_2,u,\xi).$} 
  	Define $F\in C^0(E,E)$ as $F(v)={\rm proj} (u,C(\tau_1,\tau_2,v,\xi) )$. Using Lemma~\ref{extra1}, we have 
   	\begin{eqnarray} 
&&  \hskip-0.7cm	\lVert F(v_1) - F(v_2) \rVert =\nonumber \\
&&\hskip-0.7cm =\lVert {\rm proj} (u,A+\tilde a(\tau_1,\xi) +\tilde c(\tau_2,v_1,\xi)) - {\rm proj} (u,A+\tilde a(\tau_1,\xi)+\tilde c(\tau_2,v_2,\xi)) \rVert\leq\nonumber\\
  	&& \hskip-0.7cm \leq  \lVert \tilde c(\tau_2,v_1,\xi) -\tilde c(\tau_2,v_2,\xi) \rVert \leq   L_2 \lVert v_1-v_2\rVert,\label{23}
  	\end{eqnarray}
\noindent so the existence of $v=v(\tau_1,\tau_2,s_1,s_2,u,\xi)$ with the required property (\ref{reqprop}) follows by applying the contraction mapping theorem (see e.g. Rudin \cite[Theorem 9.23]{rudin}). 
  	  	  	
\vskip0.2cm  	  	  	
  	  	  	
\noindent {\bf Step 2.} {\it Continuity of $v(\tau_1,\tau_2,s_1,s_2,u,\xi).$}   	  	  	
  	  	  	  	To prove the continuity of $v$, let $ v= v(\tau_1, \tau_2,s_1,s_2,u,\xi) \text{ and } \bar{v} = v(\bar{\tau}_1,\bar{\tau}_2,\bar{s}_1, \bar{s}_2,\bar{u},\bar{\xi})$ where $s_1,s_2,\bar{s}_1,\bar{s}_2 \in [0,T]$, $\tau_1,\tau_2,\bar{\tau}_1,\bar{\tau}_2 \in [0,T]$, $\xi, \bar{\xi}\in W$ and  $u ,\bar{u} \in E$. 

\vskip.2cm

\noindent First observe that
\begin{align*}
&\|\bar{v}-v \|=\\
 &\hskip0.2cm= \lVert {\rm proj} (\bar{u},A+\tilde a(\bar{\tau}_1,\bar{\xi})+\tilde c(\bar{\tau}_2,\bar{v},\bar{\xi}))-{\rm proj} (u ,A+\tilde a(\tau_1,\xi)+\tilde c(\tau_2,v,\xi)) \rVert \\
&\hskip0.2cm \leq \lVert  {\rm proj} (\bar{u},A+\tilde a(\bar{\tau}_1,\bar{\xi})+\tilde c(\bar{\tau}_2,\bar{v},\bar{\xi}))-{\rm proj} (u ,A+\tilde a(\bar{\tau}_1,\bar{\xi})+ \tilde c(\bar{\tau}_2,\bar{v},\bar{\xi})) \rVert\\
&\hskip0.2cm +\lVert  {\rm proj} (u,A+\tilde a(\bar{\tau}_1,\bar{\xi})+ \tilde c(\bar{\tau}_2,\bar{v},\bar{\xi})) -{\rm proj} (u ,A+\tilde a(\tau_1,\xi)+\tilde c(\tau_2,v,\xi))\rVert.
\end{align*}

\noindent Since for any nonempty, closed, convex set $C\subset E$ and any vectors $\bar u,u\in E$, we have  	(see e.g. Mordukhovich-Nam \cite[Proposition 1.79]{mord})  	  	 
  	 \begin{equation}\label{pp}
 \|{\rm proj}(\bar u,C)-{\rm proj}(u,C)\|\le \|\bar u-u\|,
 \end{equation} then,  using also Lemma~\ref{extra1}, we conclude that 
\begin{eqnarray}
\|\bar{v}-v \|&\leq& \|\bar{u}-u\|+ \left\|\tilde a(\bar{\tau}_1,\bar{\xi})+ \tilde c(\bar{\tau}_2,\bar{v},\bar{\xi})-\tilde a(\tau_1,\xi)-\tilde c(\tau_2,v,\xi) \right\| \le\nonumber\\
&\leq& \|\bar{u}-u\|+\left\| \tilde a(\bar{\tau}_1,\bar{\xi})-\tilde a(\bar{\tau}_1,\xi)\right\|+{\rm var}(\bar a,[\tau_1,\overline\tau_1])+\nonumber\\
&& +\|\tilde c(\bar{\tau}_2,v,\bar\xi)-\tilde c(\tau_2,v,\xi) \|+L_1|\bar{\tau}_2-\tau_2|+L_2\|\bar{v}-v\|,\label{24}
\end{eqnarray}   
  	\noindent so that the required continuity of  $v(\tau_1,\tau_2,s_1,s_2,u,\xi)$ follows from $0\le L_2<1$. 
  	
\vskip0.2cm

\noindent {\bf Step 3.} {\it Proof of the estimate (\ref{estimate}).} Assuming that $
  	u\in C(s_1,s_2,u,\xi)$, we have follow the lines of (\ref{24}) to get
  	\begin{eqnarray*}
  	\lVert v-u \rVert & =& \lVert {\rm proj} (u,C(\tau_1,\tau_2,v,\xi) ) - u \rVert =\min\limits_{\bar v\in  C(\tau_1,\tau_2,v,\xi)} \|u-\bar v\|.
  \end{eqnarray*}
 But $C(s_1,s_2,u,\xi)=A+\tilde a(s_1,\xi)+\tilde c(s_2,u,\xi)$ and $C(\tau_1,\tau_2,v,\xi)=A+\tilde a(\tau_1,\xi)+\tilde c(\tau_2,v,\xi)$.  Therefore,	
 \begin{eqnarray}
 \min\limits_{\bar v\in  C(\tau_1,\tau_2,v,\xi)} \|u-\bar v\|&\le& \|\tilde a(s_1,\xi)+\tilde c(s_2,u,\xi)-\tilde a(\tau_1,\xi)-\tilde c(\tau_2,v,\xi)\|\le\nonumber\\
  	& \leq & {\rm var}(\bar a,[s_1,\tau_1])  +L_1 |\tau_2-s_2|+L_2 \lVert u-v\rVert,\label{min}
  	\end{eqnarray} 
  	which implies (\ref{estimate}).

\vskip0.2cm

\noindent The proof of the lemma is complete.
  \qed

\begin{rem}\label{remivan} {\it On the validity of Lemma~\ref{lem7} when $A+c(t,\xi)$ is replaced by a more a more general term $A(t,\xi)$.} \end{rem}

\begin{floatingfigure}[r]{4.5cm} 
\hskip-0.5cm\includegraphics[scale=0.9]{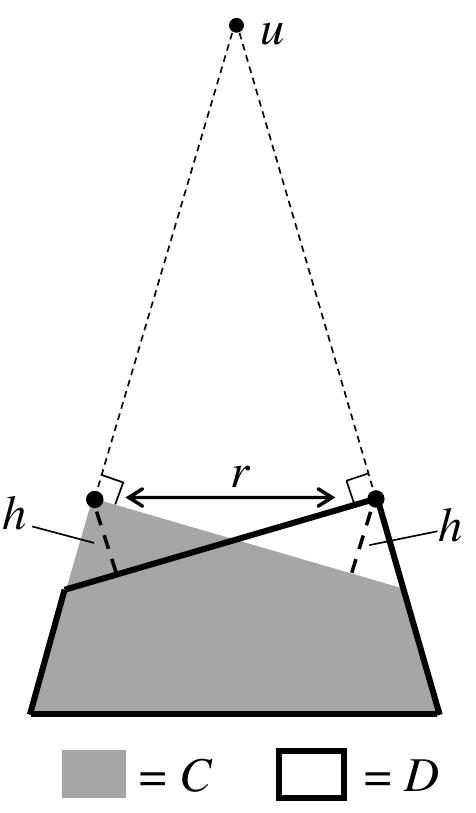} 
\caption{Illustration of  the  incorrectness of formula (\ref{contr}). Here $h=d_H(C,D),$ and  \\ $r=\|{\rm proj}(u,C)-{\rm proj}(u,D)\|$.} 
\label{contrfigure} 
\end{floatingfigure} 

\vskip-0.3cm

\noindent One can observe that estimate (\ref{23}) holds also in the case where $A+\tilde a(t,\xi)$ takes a more general form $A(t,\xi)$. Furthermore, if $d_H(A_1,A_2)$ is the Hausdorff distance between nonempty closed sets $A_1,A_2\subset E$ and $A(t,\xi)$ satisfies

\begin{equation}\label{newA1}
    d_H(A(s),A(t))\le {\rm var}(\bar a,[s,t]),
\end{equation}
then (\ref{min}) holds as well since
\begin{eqnarray*}
   \min\limits_{\bar v\in  C(\tau_1,\tau_2,v,\xi)} \|u-\bar v\|\le d_H(A(s_1,\xi)+\\
   +\tilde c(s_2,u,\xi),A(\tau_1,\xi)+\tilde c(\tau_2,v,\xi)).
\end{eqnarray*}
To summarize, the existence of  $v(\tau_1,\tau_2,s_1,s_2,u,\xi)$ (Step~1) and the estimate (\ref{estimate}) (Step~3) still hold, if $A+\tilde a(t,\xi)$ is replaced by $A(t,\xi)$ satisfying (\ref{newA1}). 

\vskip0.2cm

\noindent At the same time, the continuity of the function $v(\tau_1,\tau_2,s_1,s_2,u,\xi)$ can no longer be established when $A+\tilde a(t,\xi)$ is replaced by $A(t,\xi)$. Indeed, the core of estimate (\ref{24}) is Lemma~\ref{extra1} which doesn't allow a generalization when $C+c$ is replaced by an arbitrary set $D.$ One might be tempted to believe that the conclusion of Lemma~\ref{extra1} can be replaced by 
\begin{equation}\label{contr}
   \|{\rm proj}(u,C)-{\rm proj}(u,D)\|\le d_H(C,D),
\end{equation}
when $C+c$ is replaced by just $D$, but formula (\ref{contr}) appears to be wrong as our Fig.~\ref{contrfigure} illustrates.

\vskip0.2cm

\noindent On the other hand Monteiro Marques \cite[Proposition~4.7, p.~26]{mmbook} implies  that
\begin{eqnarray}
 && \hskip-1.8cm\|{\rm proj}(u,C)-{\rm proj}(u,D)\|\le \nonumber\\
 && \hskip-1cm\le 
 \sqrt{2({\rm dist}(u,C)+{\rm dist}(u,D))}\cdot \label{mmformula} \\
 && \cdot
 \sqrt{d_H(C,D)}, \nonumber
\end{eqnarray}
which could potentially help to obtain other versions of Lemma~\ref{lem7}, that we don't pursue in this paper.

\begin{cor}\label{corcatchup} Assume that condition I) of Theorem~\ref{thm3} holds.  Then, for any $(q,\lambda)$ satisfying (\ref{q0lam}) the implicit scheme (\ref{c1})-(\ref{c4}) is solvable iteratively from $i=0$ to $i=n-1$ and the respective iterations $x_i^n=x_i^n(q,\lambda)$ and $u_i^n=u_i^n(q,\lambda)$ are continuous in $(q,\lambda)$ on $E\times[0,1].$  Moreover,
$$
  \|u_{i+1}^n(q,\lambda)-u_i^n(q,\lambda)\|\le\dfrac{{\rm var}(\bar a,[t_i,t_{i+1}])+L_1T/n}{1-L_2},\quad i\in\overline{0,n-1},
$$
where $L_1>0$ and $L_2\in(0,1).$
\end{cor}

\noindent {\bf Proof.} Let $\xi = ((\xi_1, \xi_2, \cdots , \xi_{n+1}),\xi_{n+2})\in E^{n+1}\times\mathbb{R}$ be defined as
$$
  \xi_{i}=x_{i-1}^n,\quad i\in\overline{1,n+1} \quad , \quad \xi_{n+2} = \lambda.
$$
Therefore, the rule (\ref{c4}) defines a function $\Psi:E^{n+1} \times \mathbb{R}\to C^0([0,T],E)$ that relates $\xi\in E^{n+1} \times \mathbb{R}$ to a piecewise linear function $x_n(t)$ defined on $[0,T]$. The statement of the Corollary 1 now follows by applying Lemma~\ref{lem7} with
\begin{eqnarray*}
 && \hskip-0.6cm \tilde c(s,u,\xi)= \left(u-\int\limits_0^{s} f(\tau,\Psi(\xi)(\tau),\xi_{n+2})d\tau,\xi_{n+2}\right)+\int\limits_0^{s} f(\tau,\Psi(\xi)(\tau),\xi_{n+2})d\tau, \\
 && \hskip-0.2cm \tilde a(s,\xi) = a(s,\xi_{n+2}).
\end{eqnarray*}
The proof of the corollary is complete.
\qed

\subsection{The Poincare map associated to the catching-up algorithm}\label{poincaremapsec}

\noindent Even though we cannot ensure the existence of a Poincare map for sweeping process (\ref{sw1lam}), we can associate the following Poincare map
$$
  P^{\lambda,n}(q)=x_n(T)
$$
to the approximations $x_n$ of the catching-up algorithm (\ref{c1})-(\ref{c4}).
Corollary~\ref{corcatchup} allows to formulate the following property of the map $P^{\lambda,n}.$

\begin{cor}\label{corcontinuity} Assume that condition I) of Theorem~\ref{thm3} holds.  Consider an open bounded set $Q\subset E$. Then, for each fixed $\lambda\in[0,1]$ and $n\in\mathbb{N}$, the Poincare map $q\mapsto P^{\lambda,n}(q)$ is continuous on $\overline Q.$
\end{cor}

\subsection{Convergence of the catching-up algorithm}\label{sec64} 

\vskip0.2cm

\noindent Let $\left(u_{n}(t,q,\lambda),x_{n}(t,q,\lambda)\right)$ be the solution $(u_n(t),x_n(t))$ of the catching-up algorithm (\ref{c1})-(\ref{c4}) with the parameter $\lambda\in[0,1]$ and the initial condition $u_n(0)=x_n(0)=q$.

\begin{lem} \label{lemconverge} Assume that condition I) of Theorem~\ref{thm3} holds.  
Consider a sequence $(\lambda_n,q_n)\to (\lambda_0,q_0)$ as $n\to\infty$ of $[0,1]\times E$ satisfying (\ref{q0lam}) for each $n\in\mathbb{N}.$ Then, there exists a subsequence $\{n_k\}_{k\in\mathbb{N}}$ such that  $\{u_{n_k}(t,q_{n_k},\lambda_{n_k})\}_{k\in\mathbb{N}}$ 
converges as $k\to\infty$ uniformly in $t\in[0,T].$
\end{lem}

\noindent{\bf Proof.} {\bf Step 1.} {\it Boundedness of} $ \{u_n(t,q_n,\lambda_n)\}_{n\in \mathbb{N}}$. Let $u_i^n$, $i\in\overline{0,n},$  be the approximations given by (\ref{c1})-(\ref{c4}) with $q=q_n$ and $\lambda=\lambda_n.$
By Corollary~\ref{corcatchup}, 
  	 \begin{eqnarray*}
  		\lVert u_n(t,q_n,\lambda_n)\rVert & \leq & \|q_n\|+ \frac{1}{(1-L_2)}\left( {\rm var}(\bar a,[0, T])+L_1 T \right),
  	 \end{eqnarray*}
so the sequence $ \{u_n(t,q_n,\lambda_n)\}_{n\in \mathbb{N}} $ is  bounded uniformly on $[0,T].$

\vskip0.2cm

\noindent {\bf Step 2.} {\it Equicontinuity of }  $ \{u_n(t,q_n,\lambda_n)\}_{n\in \mathbb{N}}$. Fix $\eps>0$. Since ${\rm var}(\bar a,[s,t])\to 0$ as $|s-t|\to 0$ (see e.g. Lojasiewicz \cite[Theorem~1.3.4, p.~16]{lo}), we can choose $\delta_1>0$ such that 
\begin{equation}\label{by}
   \dfrac{{\rm var}(\bar a,[s,t])+L_1(t-s)}{1-L_2}<\dfrac{\eps}{3},\quad \mbox{for all}\ 0\le s\le t\le T\ {\rm with}\ t-s<\delta_1.
\end{equation}
Fix some $0\le s\le t\le T$ satisfying $t-s<\delta_1$ and denote by $i_s,i_t\in\overline{0,n-1}$ such indexes that
$$
   s\in[t_{i_s},t_{i_s+1}],\quad t\in[t_{i_t},t_{i_t+1}].
$$
Then we can estimate $\|u_n(t)-u_n(s)\|$ as follows:
\begin{eqnarray*}
&&   \|u_n(t)-u_n(s)\|\le\\
&& \qquad\le\|u_n(s)-u_n(t_{i_s+1})\|+\|u_n(t_{i_s+1})-u_n(t_{i_t})\|+\|u_n(t_{i_t})-u_n(t)\|\le\\
&& \qquad \le {\rm var}(u_n,[t_{i_s},t_{i_s+1}])+{\rm var}(u_n,[t_{i_s+1},t_{i_t}])+
{\rm var}(u_n,[t_{i_t},t_{i_t+1}]).
\end{eqnarray*}
The second term is smaller than $\eps/3$ by (\ref{by}) right away.  Assuming that $n\ge T/\delta_1,$ the property (\ref{by}) ensures that first and third terms are each smaller than $\eps/3$ as well. So we proved that
$$
   \|u_n(t)-u_n(s)\|< \eps,\quad \mbox{for all}\ 0\le s\le t\le T\ {\rm with}\ t-s<\delta_1,\ {\rm and}\ n\ge T/\delta_1.
$$
Since there is only a finite number of $n\in\mathbb{N}$ with $n<T/\delta_1,$ we can find $\delta_2>0$ such that 
$$
   \|u_n(t)-u_n(s)\|< \eps,\quad \mbox{for all}\ 0\le s\le t\le T\ {\rm with}\ t-s<\delta_2,\ {\rm and}\ n< T/\delta_1.
$$
Letting $\delta=\min\{\delta_1,\delta_2\}$, we finally obtain
$$
   \|u_n(t)-u_n(s)\|< \eps,\quad \mbox{for all}\ 0\le s\le t\le T\ {\rm with}\ t-s<\delta,\ {\rm and}\ n\in\mathbb{N}.
$$

\noindent The conclusion of the Lemma now follows by applying the Arzela-Ascoli theorem (see e.g. Rudin \cite[Theorem 7.25]{rudin}).
	\qed

\sloppy

\begin{rem} Establishing the existence of a converging subsequence $\{x_{n_k}(t,q_{n_k},\lambda_{n_k})\}_{k\in\mathbb{N}}$ needs more work compared to what we did in the proof of Lemma~\ref{lemconverge} because the direct corollary  of (\ref{c2}) 
$$
  x_{i+1}^n-x_i^n=u_{i+1}^n-u_i^n+\int_{t_{i-1}}^{t_i}f(\tau,x_n(\tau),\lambda)d\tau
$$
doesn't imply uniform boundedness of $x_{n}(t,q_{n},\lambda_{n}),$ $n\in\mathbb{N}$, directly. 
\end{rem}

\noindent To prove the convergence of $\{x_{n_k}(t,q_{n_k},\lambda_{n_k})\}_{k\in\mathbb{N}}$ we will now extend the discrete map (\ref{c3}) to such an operator $F_n:C([0,T],E)\to C([0,T],E)$ whose fixed point is exactly 
$t\mapsto x_n(t,q_n,\lambda_n).$ The convergence of $x_{n_k}$ will then follow from the continuity of $F_n$ in $n$ at $n=\infty.$

\vskip0.2cm

\noindent Let us define $P_n:C([0,T],E)\to E^{n+1}$, $l^-:E^{n+1}\to E^{n+1}$ and $Q_n:E^{n+1}\to C([0,T],E^{n+1})$ as
\begin{eqnarray*}
P_n(x)&=&\left(x(0),x\left(\frac{T}{n}\right),...,x\left((n-1)\frac{T}{n}\right),x(T)\right),\  \ x\in C([0,T],E),\\\\
{[l^-(y)]}_1&=&0\ ,\qquad {[l^-(y)]}_i\ =\ y_{i-1}\ , \ \  i\in\overline{2,n+1}\ ,\ y\in {E}^{n+1}\ ,\\\\
   Q_n(y)(t) &=& \dfrac{t-t_{i-1}}{1/n}y_{i+1} + \dfrac{t_{i}-t}{1/n}y_{i}\ ,\ \  y\in {E}^{n+1}\ ,\ t\in [t_{i-1},t_{i}),\ i\in\overline{1,n},\\\\
   Q_n(y)(t_n) &=& y_{n+1}\ ,\ \  y\in E^{n+1},\ \mbox{since }t_n=T.
\end{eqnarray*}
For a fixed $\lambda\in[0,1]$ and a continuous function $u:[0,T]\to E$, we introduce a continuous extension of (\ref{c3}) as
\begin{equation}\label{F}
\begin{array}{rcl}
  (F_n x)(t)&=&(Q_n P_n u)(t)-\left(Q_n l^- P_n J\right)(t),\quad t\in[0,T],\\
& & \mbox{where } J(t)=\int\limits_0^t f(\tau,x(\tau),\lambda)d\tau.
\end{array}
\end{equation}
Then, for $x\in C([0,T],E)$ satisfying $x=F_nx$, one has
\begin{eqnarray*}
x(0)&=&(Q_n P_n u)(0)-\left(Q_n l^- P_n J\right)(0)=[P_nu]_1-[l^- P_n J]_1=u(0)-0,\\
x(t_1)&=&[P_nu]_2-[l^-P_nJ]_2=u(t_1)-[P_nJ]_1=u(t_1)-J(0)=u(t_1),\\
x(t_2)&=&u(t_2)-J(t_1),\\
&\ldots&\\
x(t_n)&=&u(t_n)-J(t_{n-1}).
\end{eqnarray*}
Therefore, if $u_n$ and $x_n$ are given by (\ref{c1})-(\ref{c4}), then, letting $u=u_n$ in (\ref{F}), the fixed point $x$ of $F_n$ verifies $x(t_i)=x_n(t_i),$ $i\in\overline{0,n}.$ And, since the function $t\mapsto (F_n x)(t)$ is linear on $[t_i,t_{i+1}],$ $i\in\overline{0,n-1},$ we conclude $x_n=x.$ In other words, if $u$ in (\ref{F}) is given by $u=u_n$, then $x_n$ is the unique fixed point of $F_n$.

\begin{lem}\label{lemkam1} Assume that the conditions of Lemma~\ref{lemkam} hold. Then, there exists $\alpha>0$ and  $L\in(0,1)$ such that
$$
   \|F_n(x_1)-F_n(x_2)\|^*\le L\|x_1-x_2\|^*,\quad n\in\mathbb{N},
$$ 
for any $x_1,x_2,u\in C([0,T],E)$, $\lambda\in[0,1],$ and
$$
   \|x\|^*=\max_{t\in[0,T]}e^{-\alpha t}\|x(t)\|.
$$
Moreover, for each $x,u\in C([0,T],E)$, and $\lambda\in[0,1],$ one has
$$
  \lim_{n\to\infty} \|F_n(x)-F(x)\|=0,\quad {\rm where}\ F(x)(t)=u(t)-\int_0^t f(\tau,x(\tau),\lambda)d\tau,
$$
where $\|\cdot\|$ is the max-norm on $[0,T]$ and $F$ is a contraction in the norm $\|\cdot\|^*.$
\end{lem}

\noindent {\bf Proof.} {\bf Step 1.} Using the definition of $Q_n$, $l^-,$ and $P_n$, we have
$$
   (Q_nl^-P_nJ)(t_{i-1})=[l^-P^nJ]_i=[P_nJ]_{i-1}=J(t_{i-2}),\quad i\in\overline{2,n+1}.
$$
So that $$ (F_nx)(t_i)=u(t_i) -J(t_{i-1}).$$
  Fix $i\in\overline{1,n-1}$ and choose any $t\in[t_i,t_{i+1}].$ Then,
\begin{eqnarray*}
&&  \|F_n(x_1)(t)-F_n(x_2)(t)\|\le \\
&& \qquad\le\max\left\{\|F_n(x_1)(t_i)-F_n(x_2)(t_i)\| , \|F_n(x_1)(t_{i+1})-F_n(x_2)(t_{i+1})\|\right\}=\\
&& \qquad=\max\left\{\left\|\int_0^{t_{i-1}}f(\tau,x_1(\tau),\lambda)d\tau-\int_0^{t_{i-1}}f(\tau,x_2(\tau),\lambda)d\tau\right\|,\right.\\
&& \qquad\qquad\qquad\left. \left\|\int_0^{t_{i}}f(\tau,x_1(\tau),\lambda)d\tau-\int_0^{t_{i}}f(\tau,x_2(\tau),\lambda)d\tau\right\|\right\}\le\\
&& \qquad\le \bar L\int_0^{t_{i}}\|x_1(\tau)-x_2(\tau)\|d\tau\le \bar L\int_0^{t_{i}}e^{\alpha \tau}\|x_1-x_2\|^*d\tau,
\end{eqnarray*}
where $\bar{L}>0$ is the global Lipschitz constant of $x\mapsto f(t,x,\lambda)$ and $\alpha>0$ is an arbitrary constant. Therefore,
$$  e^{-\alpha t}\|F_n(x_1)(t)-F_n(x_2)(t)\| \le \dfrac{\bar L}{\alpha}\left(e^{\alpha(t_{i}-t)}-e^{-\alpha t}\right)\|x_1-x_2\|^*\le \dfrac{\bar L}{\alpha}\|x_1-x_2\|^*,
$$
which holds for any $t\in[0,T].$ The case $t\in[0,t_1]$ can be considered along the same lines. This proves the contraction part of the lemma.

\vskip0.2cm

\noindent {\bf Step 2.} To prove the convergence part, fix $i\in\overline{1,n-1}$ again and consider $t\in[t_i,t_{i+1}].$  
Since $(Q_nP_nu)(t_i)=u(t_i),$ we have
\begin{eqnarray*}
  \|(Q_nP_nu)(t)-u(t)\| &\le& \|(Q_nP_nu)(t)-(Q_nP_nu)(t_i)\|+\|u(t)-u(t_i)\|\le\\
 &\le&\|u(t_{i+1})-u(t_i)\|+\|u(t)-u(t_i)\|,
\end{eqnarray*}
so that the convergence of $(Q_nP_nu)(t)$ to $u(t)$ as $n\to\infty$ follows from continuity of $u.$ The convergence of $(Q_nl^-P_nJ)(t)$ follows same lines. Indeed, since $(Q_nl^-P_nJ)(t_{i+1})=J(t_i),$ one has
\begin{eqnarray*}
 &&  \|(Q_nl^-P_nJ)(t)-J(t)\|\le\\
&&\qquad \le\|(Q_nl^-P_nJ)(t)-(Q_nl^-P_nJ)(t_{i+1})\|+\|J(t)-J(t_i)\|\le\\
&&\qquad\le \|J(t_{i-1})-J(t_i)\|+\|J(t)-J(t_i)\|
\end{eqnarray*}
and the convergence of $(Q_nl^-P_nJ)(t)$ to $J(t)$ as $n\to\infty$ follows from continuity of $J(t).$

\vskip0.2cm

\noindent The proof of the lemma is complete.\qed

\begin{cor}\label{corkam} 
 Assume that condition I) of Theorem~\ref{thm3} holds.  Let  $\{n_k\}_{k\in\mathbb{N}}$ be the subsequence given by Lemma~\ref{lemconverge} (which ensures the convergence of 
   $\{u_{n_k}(t,q_{n_k},\lambda_{n_k})\}_{k\in\mathbb{N}}$).
Consider the limit
$$
  u(t)=\lim\limits_{k\to\infty}u_{n_k}(t,q_{n_k},\lambda_{n_k}).
$$ 
Let $x(t)$ be the solution of the respective equation (\ref{xu}) (which exists according to Lemma~\ref{lemkam}). Then $\{x_{n_k}(t,q_{n_k},\lambda_{n_k})\}_{k\in\mathbb{N}}$ converges uniformly in $t\in[0,T]$, and
\begin{equation}\label{xlimit}
x(t)=\lim\limits_{k\to\infty}x_{n_k}(t,q_{n_k},\lambda_{n_k}).
\end{equation}
\end{cor}

\noindent {\bf Proof.} The conclusion follows from the inequality
\begin{eqnarray*}
   \|x-x_n\|^*&\hskip-0.1cm=&\hskip-0.1cm\|F(x)-F_n(x_n)\|^*\le\|F(x)-F_n(x)\|^*+\|F_n(x)-F_n(x_n)\|^*\le\\
&& \le \|F(x)-F_n(x)\|^*+L\|x-x_n\|^*,
\end{eqnarray*}
where $L\in(0,1)$ is given by Lemma~\ref{lemkam1}.\qed

\subsection{Verifying that the limit of the catching-up algorithm is indeed a solution}\label{sec65}

\vskip0.2cm
\begin{thm}\label{thmexistence} Let the conditions of Corollary~\ref{corkam} hold and let $u(t)$ and $x(t)$ be as given by this corollary. Then, $u(t)$ is a  solution of sweeping process (\ref{sw3})
with the parameters
$
x(t),$ $\lambda=
\lim\limits_{k\to\infty}\lambda_{n_k}
$,
and the initial condition
$ 
u(0)=\lim\limits_{k\to\infty}q_{n_k}.
$ Accordingly, by Definition~\ref{df1}, $x(t)$ is a solution of perturbed sweeping process (\ref{sw1lam}).  
\end{thm}

\noindent {\bf Proof.} Let $\phi(t)$, $t\in[0,T],$ be an arbitrary continuous selector of the moving set of (\ref{sw3}), i.e.
$$
  \phi(t)\in A+a(t,\lambda)+c(x(t),\lambda)+\int_0^t f(\tau,x(\tau),\lambda)d\tau,\quad t\in[0,T].
$$
According to Monteiro Marques \cite[p.~15-16]{26} (see also Valadier \cite[Proposition~6]{valadier}) it is sufficient to prove that
\begin{equation}\label{required}
   \int_s^t\left<\phi(\tau),du(\tau)\right>\ge \dfrac{1}{2}\left(\|u(t)\|^2-\|u(s)\|^2\right),\quad 0\le s\le t\le T,
\end{equation}
which we now establish using the ideas of Kunze and Monteiro Marques \cite{kunze2}.

\vskip0.2cm

\noindent Without loss of generality we will assume that $\{n_k\}_{k\in\mathbb{N}}=\mathbb{N},$ and replace $n_k$, $k\in\mathbb{N}$, by $n$, $n\in\mathbb{N}$ in the formulation of the theorem.
Fix $t>0$ and select $i\in\overline{0,n-1}$ such that $t\in[t_i,t_{i+1}].$ Introduce $\hat c_n(t)$ as
$$
  \hat c_n(t)={\rm proj}\left(\phi(t), A+a(t_{i+1},\lambda_n)+c(x^n_{i+1},\lambda_n)+\int_0^{t_i} f(\tau,x_n(\tau),\lambda_n)d\tau\right).
$$
Then, by (\ref{c1}) and by convexity of $A$, we have (see e.g. Kunze and Monteiro Marques \cite[formula (4)]{kunze2})
$$
   \left<u_n(t_{i+1})-u_n(t_i),u_n(t_{i+1})-\hat c_n(t)\right>\le 0, \quad t\in[t_i,t_{i+1}],
$$
from where
\begin{eqnarray*}
&&   \left<u_n(t_{i+1})-u_n(t_i),u_n(t)-\hat c_n(t)\right>\le\\
&&\qquad \le \left<u_n(t_{i+1})-u_n(t_i),u_n(t)-u_n(t_{i+1})\right>\le \lVert u_n(t_{i+1})-u_n(t)\rVert ^2,
\end{eqnarray*}
or 
\begin{eqnarray*}
&& \hskip-0.7cm \left<u_n(t_{i+1})-u_n(t_i),\hat c_n(t)\right>\ge -\lVert u_n(t_{i+1})-u_n(t_i)\rVert ^2+\left<u_n(t_{i+1})-u_n(t_i),u_n(t)\right>,
\end{eqnarray*}
for any $t\in[t_i,t_{i+1}].$ Using the linearity of $u_n$ on $[t_i,t_{i+1}]$, we conclude 
\begin{eqnarray*}
&&\hskip-0.7cm\left<\hat c_n(t),u_n(\bar t_{i+1})-u_n(\bar t_i)\right>\ge \\
&&\ge \left<u_n(t),u_n(\bar t_{i+1})-u_n(\bar t_i)\right>-\left<u_n(\bar t_{i+1})-u_n(\bar t_i)) , (u_n(t_{i+1})-u_n(t_i)) \right>,
\end{eqnarray*}
for any $t_i\le \bar t_i\le t\le\bar t_{i+1}\le t_{i+1}.$
Therefore, denoting $\tau_{j,k}=\bar t_i+\left(j+\frac{1}{2}\right)\frac{\bar t_{i+1}-\bar t_i}{k}$ for $j\in \{0,1,\cdots,k-1\}$,  one has (same approach is used e.g. in part (ii) of the proof of Monteiro Marques \cite[Theorem~2.1, p.~12, second formula from below]{mmbook})
\begin{eqnarray*}
 &&\hskip-0.7cm  \int_{\bar t_i}^{\bar t_{i+1}} \left<\hat c_n(\tau),du_n(d\tau)\right>=\\
&&\hskip-0.7cm =\lim_{k\to\infty}\sum_{j=0}^{k-1}\left<\hat c_n\left(\tau_{j,k}\right),u_n\left(\bar t_i+(j+1)\frac{\bar t_{i+1}-\bar t_i}{k}\right)-u_n\left(\bar t_i+j\frac{\bar t_{i+1}-\bar t_i}{k}\right)\right>\ge\\
&&\hskip-0.7cm\ge\lim_{k\to\infty}\sum_{j=0}^{k-1} \left<u_n\left(\tau_{j,k}\right),u_n\left(\bar t_i+(j+1)\frac{\bar t_{i+1}-\bar t_i}{k}\right)-u_n\left(\bar t_i+j\frac{\bar t_{i+1}-\bar t_i}{k}\right)\right>+R_n,
\end{eqnarray*}
where the reminder ${R_n}$ is given by
 \begin{eqnarray*} && R_n  =-\lim_{k\to\infty}\sum_{j=0}^{k-1}\left<u_n\left(\bar t_i+(j+1)\frac{\bar t_{i+1}-\bar t_i}{k}\right)-\right. \\
&& \qquad\qquad\qquad\qquad -\left.u_n\left(\bar t_i+j\frac{\bar t_{i+1}-\bar t_i}{k}\right),u_n(t_{i+1})-u_n(t_i)\right> =\\
&&\quad\ \ =-\left< u_n(\bar t_{i+1})-u_n(\bar t_i), u_n(t_{i+1})-u_n(t_i)\right>.\end{eqnarray*}

\noindent Therefore,
\begin{eqnarray*}
&& \int_{\bar t_i}^{\bar t_{i+1}} \left<\hat c_n(\tau),du_n(d\tau)\right>=\\
&& =\lim_{k\to\infty}\sum_{j=0}^{k-1}\left<u_n(\tau_{j,k}),u_n'(\tau_{j,k})\dfrac{\bar t_{i+1}-\bar t_i}{k}\right>+R_n=\\
&& =\dfrac{1}{2}\lim_{k\to\infty}\sum_{j=0}^{k-1}
\left.\left(\dfrac{d}{d\tau}\|u_n(\tau)\|^2\right)\right|_{\tau=\tau_{j,k}}\cdot\dfrac{\bar t_{i+1}-\bar t_i}{k}+R_n=\\
&& =\dfrac{1}{2}\int_{\bar t_i}^{\bar t_{i+1}}\dfrac{d}{d\tau}\|u_n(\tau)\|^2d\tau+R_n=\dfrac{1}{2}\left(\|u_n(\bar t_{i+1})\|^2-\|u_n(\bar t_i)\|^2\right)+R_n.
\end{eqnarray*}
 	
\noindent This result can now be used to estimate the required integral (\ref{required}) as follows
\begin{equation}\label{Rn}
  \int_{s}^{t} \left<\hat c_n(\tau),du_n(d\tau)\right>\ge 
         \dfrac{1}{2}\left(\|u_n(t)\|^2-\|u_n(s)\|^2\right)+R,
\end{equation}
where
$$
   |R|\le {\rm var}(u_n,[s,t])\cdot \max_{i\in\overline{0,n-1}}\|u_n(t_{i+1})-u_n(t_i)\|.
$$
But according to Corollary~\ref{corcatchup}, 
$$
   {\rm var}(u_n,[s,t])\le \dfrac{{\rm var}(\bar a,[s,t])}{1-L_2}+\dfrac{L_1|t-s|}{1-L_2}.
$$
Therefore, the desired statement (\ref{required}) follows from (\ref{Rn}) by passing to the limit as $n\to\infty$ (the passage to the limit is valid e.g. by Monteiro Marques \cite[Theorem~2.1(ii)-(iii)]{mmbook} combined with formula (26) of p.~7 of the same book). 

\vskip0.2cm

\noindent The proof of the theorem is complete. \qed

\subsection{Proof of Theorem~\ref{thm1} (sweeping process without a parameter)}

\noindent 
Theorem~\ref{thm1} is a direct consequence of Theorem~\ref{thmexistence}. One just view sweeping process (\ref{sw1}) as sweeping process (\ref{sw1lam}) with $\lambda=0$. 

\begin{rem}\label{remgeneral}
Using Remark~\ref{remivan}, Theorem~\ref{thm1} can be directly extended to sweeping processes of the form
\begin{equation}\label{swextend}
   -dx\in N_{A(t,\lambda)+c(x,\lambda)}(x)+f(t,x,\lambda)dt,\quad  x\in E, \ \lambda\in\mathbb{R},
\end{equation}
where $A$ is a set-valued function with nonempty closed convex bounded values that satisfies the property
\begin{equation}\label{A1_}
\begin{array}{l}
   d_H(A(s,\lambda),A(t,\lambda))\le {\rm var}(\bar a,[s,t]),\quad \lambda\in[0,1],\\
\mbox {where }\bar a:[0,T]\to\mathbb{R}\mbox{ is a BV continuous function}.
\end{array}
\end{equation}

\end{rem}

\vskip0.2cm 


\subsection{Proofs of Theorems~\ref{thm3} and \ref{thm4} (sweeping process with a parameter)}

\vskip0.2cm

\noindent{\bf Proof of Theorem~\ref{thm3}.} {\bf Step 1.} First we prove that there exists $N>0$ such that $d(I-P^{\lambda,n}\circ V^{\lambda,n},Q)$ is defined for $n\ge N$ and $\lambda\in[0,\lambda_1]$. Assuming the contrary, we get a sequence $n_k\to\infty$, $\lambda_k\to\lambda_0\in[0,\lambda_1]$, and a converging sequence $\{q_k\}_{k\in\mathbb{N}}\subset\partial Q$ such that 
\begin{equation}\label{step1}
   P^{\lambda_k,n_k}\circ V^{\lambda_k,n_k}(q_k)=q_k,\quad k\in\mathbb{N}.
\end{equation}
Applying Lemma~\ref{lemconverge}, Corollary~\ref{corkam}, and Theorem~\ref{thmexistence} we conclude that $q_0=\lim_{k\to\infty}q_k\in\partial Q$ is the initial condition of the $T$-periodic solution (\ref{xlimit}) of sweeping process (\ref{sw1lam}) with $\lambda=\lambda_0$, which contradicts conditions III) of Theorem~\ref{thm3}.

\vskip0.2cm

\noindent The conclusion of Step 1, in particular, implies that
$$
   d(I-P^{\lambda,n}\circ V^{\lambda,n},Q)=d(I-P^{0,n}\circ V^{0,n},Q),\quad n\ge N,\ \lambda\in[0,\lambda_1].
$$

\noindent {\bf Step 2.} Here we use assumption II (uniqueness) of Theorem~\ref{thm3} to conclude that 
$$
   P^{0,n}\circ V^{0,n}(q)\to P^0\circ V^{0}(q),\quad {\rm as}\ n\to\infty,
$$
uniformly with respect to $q\in\overline{Q}.$ Thus, we can  diminish $N>0$ in such a way that $d(I-P^{0,n}\circ V^{0,n},Q)=d(I-P^0\circ V^{0},Q)$, $n\ge N,$ which gives
$$
   d(I-P^{\lambda,n}\circ V^{\lambda,n},Q)\not=0,\quad n\ge N, \ \lambda\in[0,\lambda_1].
$$
Therefore, for each $\lambda\in[0,\lambda_1]$ there exists $q_n\in Q$ such that the approximations   $\{x_n(\cdot,q_n,\lambda)\}_{n\ge N}$ are $T$-periodic, so this sequence has a convergent subsequence which converges to a $T$-periodic solution of (\ref{sw1lam}) with initial condition $q = \lim_{n_k\to \infty}q_{n_k}$ as $n\to\infty$ according to Corollary~\ref{corkam}. 

\vskip0.2cm

\noindent The proof of the theorem is complete.\qed

\vskip0.2cm

\noindent The proof of Theorem~\ref{thm4} follows the lines of the proof of Theorem~\ref{thm3}. The only difference is in the beginning of  Step~1, which now proves the existence of both $N>0$ and $\lambda_1\in(0,1]$ such that $d(I-P^{\lambda,n}\circ V^{\lambda,n},Q)$ is defined for $n\ge N$ and $\lambda\in[0,\lambda_1]$. Assuming the contrary, we get a sequence $n_k\to\infty$, $\lambda_k\to 0\in[0,1]$, and a converging sequence $\{q_k\}_{k\in\mathbb{N}}\subset\partial Q$ such that (\ref{step1}) holds, that leads to the existence of a $T$-periodic solution to sweeping process (\ref{sw1lam}) with $\lambda=0,$ contradicting condition II) of Theorem~\ref{thm3}. The rest of the proof of Theorem~\ref{thm4} follows the proof of Theorem~\ref{thm3} just literally.

\section{Proof of the theorem on the global existence of periodic solutions}

\noindent To prove Theorem~\ref{thmexist} we will use the following well-known result (see e.g. Krasnoselskii-Zabreiko \cite[Theorem~6.2]{kz}):

\begin{thm}\label{thmkra} Let $\bar P:E\to E$ be a continuous map and let  $Q\subset E$ be an open bounded convex set. If $\bar P(Q)\subset\overline{Q}$ and if $\bar P$ doesn't have fixed points on $\partial Q$,  then
$$
   d(I-\bar P,Q)=1.
$$
\end{thm}

\noindent{\bf Proof of Theorem~\ref{thmexist}.} Let $\Omega_1$ be the  1-neighborhood of $\Omega.$ Since $\Omega$ is convex, then $\Omega_1$ is convex as well. 
 We will view sweeping process (\ref{sw1}) as sweeping process (\ref{sw1lam}) with $\lambda=0$. 
So we consider the map 
$$
   \overline{P}{}^{0,n}(x)=P^{0,n}(V(x)),
$$
where $P^{0,n}$ is as introduced in Section~\ref{poincaremapsec} and $V$ is as introduced in Section~\ref{3}. We claim that 
\begin{equation}\label{abcd}
 \overline P^{0,n}(\Omega_1)\subset\Omega,\quad\mbox{for all }n\in\mathbb{N}.
\end{equation}
We have $V(x)\in\Omega$ by the definition of the map $V.$ Then, according to the catching-up scheme~(\ref{c1})-(\ref{c4}), we have that
$$
   x_{i+1}^n\in A+a(t_{i+1},0)+c(x_{i+1}^n),\quad{\rm i.e.}\ x_{i+1}^n\in\Omega_{t_{i+1}},\ i\in\overline{0,n-1},
$$
and so $x_n(T)\in\Omega_T$, which implies (\ref{abcd}).

\vskip0.2cm

\noindent Using the continuity of $P^{0,n}$ (Corollary~\ref{corcontinuity}) and $V$ (Lemma~\ref{lem7}) along with Theorem~\ref{thmkra}, we get the existence of $q_n\in\Omega$ such that
$
   \overline P{}^{0,n}(q_n)=q_n,
$
which implies 
$$
   P{}^{0,n}(q_n)=q_n,\quad n\in\mathbb{N},
$$
because $V(q_n)\in\Omega.$ In other words, we have $x_n(T,q_n,0)=x_n(0,q_n,0)$ for all $n\in\mathbb{N}.$ Now, Theorem~\ref{thmexistence} applied with $\lambda_n=0$, implies the existence of a convergent subsequence $\{x_{n_k}(t,q_{n_k},0)\}$ whose limit $x(t)$ is solution of (\ref{sw1}) with the required $T$-periodicity property (\ref{Tper}).  The proof is complete. \qed

\section{Existence of periodic solutions in the neighborhood of a boundary equilibrium (the theorem and its proof)} \label{7}

\noindent This section uses the following extension of Theorem~\ref{thmkra} (see e.g. Krasnoselskii-Zabreiko \cite[Theorem~31.1]{kz}):

\begin{thm}\label{thmkra1} Let $\bar P:E\to E$ be a continuous map and let  $Q\subset E$ be an open bounded set. If $(\bar P)^m$ maps $Q$ strictly into itself for all $m\in\mathbb{N}$ sufficiently large, then
$$
   d(I-\bar P,Q)=1.
$$
\end{thm}

\noindent The main assumption of this section is that sweeping processes (\ref{sw1lam}) reduces to
\begin{equation}\label{swauto}
   -\dot x\in N_{A}(x)+f_0(x), \qquad x\in E,
\end{equation}
when $\lambda=0$ and that (\ref{swauto}) posses a switched equilibrium on the boundary $\partial A$ (as was earlier introduced in Kamenskii-Makarenkov \cite{kamenskii} in 2d). To introduce the definition of a switched boundary equilibrium $x_0\in\partial A$, we assume that in some neighborhood $Q\subset \mathbb{R}^n$ of $x_0$ the boundary $\partial A$ is smooth and can be described as
$$
   \partial A\cap Q=\{x\in Q: H(x)=0\}, \quad {\rm where}\ H\in C^1(\mathbb{R}^n,\mathbb{R}).
$$ 

\begin{df}\label{se} A point $x_0\in\partial A$ is a switched boundary equilibrium of sweeping process (\ref{swauto}), if 
$$
   H(x)>0,\quad\mbox{for all}\ x\in Q\backslash A,
$$
and
$$
  H'(x_0)=\alpha f(x_0)\quad \mbox{for some}\ \alpha<0.
$$
\end{df}
 As the definition says, $x_0$ is not an equilibrium of $f$, however the next two lemmas imply that the solution of (\ref{swauto}) with the initial condition at $x_0$ don't leave $x_0.$
 
 \vskip0.2cm
 
 \noindent If $x_0$ is a switched equilibrium, then $Q$ can be considered so small that
\begin{equation}\label{sliding}
    \left<f(x),H'(x)\right><0,\quad\mbox{for all}\ x\in \partial A\cap Q.
\end{equation}

\noindent The next lemma claims that $\partial A\cap Q$ is a sliding region for sweeping process (\ref{swauto}).

\begin{lem}\label{lemslide} Let $x_0\in\partial A$ be a switched equilibrium of (\ref{swauto}) and let $Q\subset E$ be such a neighborhood of $x_0$ that (\ref{sliding}) holds. Consider a solution $x$ of (\ref{swauto}) with an initial condition $x_0\in \partial A\cap Q.$ Let $t_1>0$ be such that $x(t)\in Q$ for all $t\in[0,t_1]$. Then $x(t)\in\partial A$ for all $t\in[0,t_1].$ 
\end{lem}
\noindent {\bf Proof.}
	Let us assume, by contradiction, that there exists $t_{escape} \in [0,t_1]$ where $x(t)$ escapes from $\partial A$, i.e.
\begin{equation*}
t_{escape} = \max\{t_0\ge 0:x(t)\in Q, H(x(t))=0,\ t\in[0,t_0]\}<t_1.
\end{equation*}
\noindent By the definition of $t_{escape},$ for any $\delta>0$ there exist $t_\delta\in[t_{escape},t_{escape}+\delta]$ such that $H(x(t))<0$ for each $t\in (t_{escape},t_\delta]$.
Since, the solution $x(t)$ satisfies $\dot{x}(t) = - f_0(x(t))$ on $(t_{escape},t_\delta]$, by the Mean-Value Theorem
$$
H(x(t_\delta))-H(x(t_{escape}) = -H'(x(t_{\delta}^*))
f_0(x(t_{\delta}^*)) (t_\delta-t_{escape}),
$$ 
for some $t_{\delta}^*\in (t_{escape},t_\delta)$.
This yields
\begin{equation*}
 H'(x(t_{escape})) f_0(x(t_{escape}))  \geq 0,
\end{equation*}
as $\delta\to 0,$ 
contradicting (\ref{sliding}). 

\vskip0.2cm

\noindent The proof of the lemma is complete.
\qed

\vskip0.2cm

\noindent As it happens in the theory of Filippov systems (see \cite{fil}), the dynamics of (\ref{swauto}) in the sliding region is described by a smooth differential equation. Indeed, let us introduce the differential equation
\begin{equation}\label{slidingeq}
 \begin{array}{l}
   -\dot x=\bar f(x),\\
{\rm where}\ \bar f(x)=  f_0(x)-\pi_{H'(x)}(f_0(x)) \ {\rm and}\   \pi_{L}(\xi)=\dfrac{1}{\|L\|^2}\left<\xi,L\right>L.
  \end{array}
\end{equation}
Next lemma says that (\ref{slidingeq}) is the equation of sliding motion for sweeping process (\ref{swauto}) in the neighborhood of switched equilibrium $x_0\in\partial A_0.$

\begin{lem}\label{leminv} Let the conditions of Lemma~\ref{lemslide} hold and let $x(t)$ be the sliding solution $x(t)$, $t\in[0,t_1],$ of sweeping process (\ref{swauto}) as introduced in Lemma~\ref{lemslide}. Then $x(t)$ is a solution of (\ref{slidingeq}) on $[0,t_1].$
\end{lem}

\noindent {\bf Proof.} Fix $t\in[0,t_1]$ such that $\dot x(t)$ exists. Then, from (\ref{swauto}), 
$$
   -\dot x(t)=\alpha H'(x(t))+f_0(x(t)),\quad\mbox{with some}\ \alpha>0,
$$
or
\begin{equation}\label{or}
   \alpha H'(x(t))=-\pi_{H'(x(t))}(f_0(x(t)))+\left[-f_0(x(t))+\pi_{H'(x(t))}(f_0(x(t)))\right]-\dot x(t).
\end{equation}
From the definition of $\pi_L(\xi)$ we have
$$
  \left<-f_0(x(t))+\pi_{H'(x(t))}(f_0(x(t))),H'(x(t))\right>=0.
$$
On the other hand, from Lemma~\ref{lemslide},
$$
   \left<\dot x(t),H'(x(t))\right>=0.
$$
Therefore, taking the scalar product of (\ref{or}) with $H'(x(t))$, we get
$$
   \alpha=-\dfrac{1}{\|H'(x(t))\|^2}\left<f_0(x(t)),H'(x(t))\right>,
$$ 
which completes the proof.\qed

\vskip0.2cm

\noindent Lemma~\ref{leminv} implies that the boundary $\partial A$ is an invariant manifold for the differential equation (\ref{slidingeq}). The definition (\ref{slidingeq}) reduces the dimension of the image of $f_0$ by 1. Therefore, the image of the map $\bar f$ acts to a space of dimension $dimE-1$, which implies that one eigenvalue of the Jacobian $\bar f'(x_0)$ is always zero.

\vskip0.2cm

\noindent We now offer an asymptotic stability result which can be of independent interest in applications of perturbed sweeping processes.

\begin{thm}\label{stabslide} Let $x_0\in\partial A$ be a switched equilibrium of (\ref{swauto}). If real parts of $dimE-1$  eigenvalues of the Jacobian $\bar f'(x_0)$ are negative, then $x_0$ is a uniformly asymptotically stable point of sweeping process (\ref{swauto}).
\end{thm}
\noindent {\bf Proof.} {\bf Step 1.} {\it Convergence to $\partial A.$} Let $B_r(x_0)$ be a ball of radius $r$ centered at $x_0.$ Let us show that there exists $r>0$ such that for any $\xi\in B_r(x_0)\cap A,$ 
   the solution $t\mapsto X(t,\xi)$ of
\begin{equation}\label{aa}
   \dot x=-f_0(x)
\end{equation}
with the initial condition $X(0,\xi)=\xi$ reaches $\partial A$ at time some time $\tau(\xi)>0$. 
   The proof will be through the Implicit Function Theorem applied to
$$
   F(t,x)=H(X(t,x)).
$$
We have $F(0,x_0)=0$ and $F_t(0,x_0)=-H'(x_0)f_0(x_0)\not=0$ by the definition of switched equilibrium. Therefore, Implicit Function Theorem (see e.g. Rudin \cite[Theorem 9.28]{rudin}) ensures the existence of $\xi\to \tau(\xi)$ defined and continuous on a sufficiently small ball $B_r(x_0)$ and such that $\tau(x_0)=0.$

\vskip0.2cm

\noindent It remains to show that $\tau(\xi)>0$ for all $\xi\in B_r(x_0)\cap A.$ Since, according to the definition of switched equilibrium, $H'(x_0)^T$ is a normal to $A$ pointing outwards to $A$, it is sufficient to prove that $\tau(\xi)>0$ for $\xi=x_0-\lambda H'(x_0)^T$ with all $\lambda>0$ sufficiently small. So we introduce a scalar function
$$
   G(\lambda)=\tau(x_0-\lambda H'(x_0)^T)
$$
and want to prove that $G'(0)>0.$ Using the formula for the derivative of the implicit function (see \cite[Theorem 9.28]{rudin}) 
$$
   \tau'(x_0)=-(H'(x_0)f_0(x_0))^{-1}H'(x_0)
$$
and so $$G'(0)=-(H'(x_0)f_0(x_0))^{-1}H'(x_0)(-H'(x_0)^T)=H'(x_0)f_0(x_0)\|H'(x_0)\|^2,$$
which is indeed positive according to Definition~\ref{se}.

\vskip0.2cm

\noindent  Finally, let us fix $\xi\in B_r(x_0)\cap A$ and let $x(t)$ be the solution of (\ref{swauto}) with the initial condition $x(0)=\xi.$ 
Since the conclusion of the Implicit Function Theorem comes with uniqueness, we have that $X(t,\xi)\not\in\partial A,$ $t\in[0,\tau(\xi)).$  Therefore, $X(t,\xi)=x(t)$, for any $t\in[0,\tau(\xi)),$ which implies that $\lim_{t\to \tau(\xi)}X(t,\xi)=\lim_{t\to \tau(\xi)}x(t)$ and so $x(\tau(\xi))\in\partial A.$ 

\vskip0.2cm

\noindent {\bf Step 2.} {\it Convergence along $\partial A.$} Lemmas \ref{lemslide} and 
\ref{leminv} combined with the negativeness of real parts of $dimE-1$ eigenvalues of $\bar f'(x_0)$ imply that there exists an neighborhood $x_0\in Q\subset E$ such that any solution of (\ref{swauto}) with the initial condition $x(0)\in Q\cap\partial A$ converges to $x_0$ along $\partial A$ as $t\to\infty$ and the convergence is uniform with respect to the initial condition.

\vskip0.2cm

\noindent Making now $r>0$ in Step~1 so small that $\cup_{\xi\in B_r(x_0)} X(\tau(\xi),\xi)\in Q$ (which is possible by continuity of $\xi\to\tau(\xi)$), we combine Step 1 and Step 2 to conclude that any solution of (\ref{swauto}) with $x(0)\in B_r(x_0)$ approaches $x_0$ as $t\to\infty.$

\vskip0.2cm

\noindent The proof of the theorem is complete.
\qed

\vskip0.2cm

\noindent We are now in the position to combine theorems~\ref{thm4}, \ref{thmkra}, and \ref{stabslide} when the following condition holds for (\ref{sw1lam}) at $\lambda=0:$
\begin{equation}\label{case1}
   a(t,0)\equiv 0,\quad c(x,0)\equiv 0,\quad f(t,x,0)\equiv f_0(x) \ {\rm with} \ f_0\in C^1(E,E).
\end{equation}

\begin{thm} \label{thm16} Assume that condition I) of Theorem~\ref{thm3} holds. Assume, that for $\lambda=0$ sweeping process (\ref{sw1lam}) is smooth autonomous, i.e. satisfies (\ref{case1}). If real parts of $n-1$ eigenvalues of $\bar f'(x_0)$ are negative for some switched equilibrium $x_0\in\partial A,$ then there exists $\lambda_1>0$ such that for all $\lambda\in(0,\lambda_1]$ sweeping process (\ref{sw1lam}) admits a periodic solution $x_\lambda(t)\to x_0$ as $\lambda\to 0.$
\end{thm}

\noindent {\bf Proof.} Let $\bar P(x) = P^0(V^0(x)).$ By Theorem~\ref{stabslide}, there exists an open bounded set $x_0\in Q\subset E$ such that  $(\bar P)^m$ maps $Q$ strictly into itself for all $m\in\mathbb{N}$ sufficiently large. Therefore, Theorem~\ref{thmkra1} ensures that condition II) of Theorem~\ref{thm3} holds, so Theorem~\ref{thm4} applies.\qed

\vskip0.2cm

\noindent Similar to Theorem~\ref{thm16} results have been obtained for ordinary differential equations by  Berstein-Halanai \cite{berstein}  and Cronin \cite{cronin}.

\section{Conclusions} \noindent By extending the implicit catching-up scheme of Kunze and Monteiro Marques \cite{kunze2} to perturbed sweeping processes, we proved solvability of BV-continuous state-dependent sweeping processes with Lipschitz dependence on the state. We further used topological degree arguments to establish the existence of periodic solutions to sweeping processes of this type. The analysis is carried out for the simplest possible moving set $C(t)=A+a(t)+c(x)$ throughout the entire paper, that allowed us to focus on the development of core mathematical ideas rather than on its possible generalizations. We explain in Remarks~\ref{remivan} and \ref{remgeneral} how the existence result (Theorem~\ref{thm1}) immediately extends to the moving set of the form  $C(t)=A(t)+c(x)$. At the same time, Remark~\ref{remivan} shows that  our method of proof of continuity of approximations $x_n(t)$ on the initial condition fails for moving sets of the form $C(t)=A(t)+c(x)$. Since continuous dependence of $x_n(t)$ on the initial condition is the main ingredient in our proof of the existence of periodic solutions, the respective main theorems (Theorem~\ref{thmexist} and Theorem~\ref{thm16}) do not readily extend even to the moving set of the form  $C(t)=A(t)+c(x)$. We don't know whether or not an alternative approach (e.g. formula (\ref{mmformula}) quoted from 
\cite[Proposition~4.7, p.~26]{mmbook})  can deal with any more general state-dependent moving constraints. 

\vskip0.2cm

\noindent The existence of a $T$-periodic solutions to a sweeping process with $T$-periodic right-hand-sides and convex moving set would be an immediate result when uniqueness and continuous dependence of solutions on initial conditions holds. The difficulty we overcame when proving the existence of periodic solutions comes from the fact that uniqueness and continuous dependence on initial conditions of solutions of BV-continuous state-dependent sweeping processes is still an open question even when the dependence on the state is Lipschitz continuous (for state-independent sweeping processes uniqueness and continuous dependence  is established e.g in Castaing and Monteiro Marques \cite{castaing} and Adly et al \cite{adlynew}).

\vskip0.2cm

\noindent The second part of the paper concerns sweeping processes with a parameter $\lambda,$ for which we developed a topological degree based continuation principle. As an application of the continuation principle, we proved the occurrence of periodic solutions at a specific location being a neighborhood of a switched boundary equilibrium. Specifically, we assumed that for $\lambda=0$, the sweeping process is autonomous and admits an asymptotically stable switched boundary equilibrium $x_0.$ We then proved the occurrence of $T$-periodic solutions from $x_0$ when the parameter $\lambda$ increases and the sweeping process becomes nonautonomous (and $T$-periodic). The condition for asymptotic stability of $x_0$ can be replaced by assuming that the topological index of $x_0$ is different from $0.$ Such a condition can be also expressed in terms of the eigenvalues of the linearization $\bar f(x_0)$ of sliding differential equation (\ref{slidingeq}), see e.g. Krasnoselskii-Zabreiko \cite[Theorem~6.1]{kz} and \cite[Theorem~7.4]{kz} (which will be required to account for the vector field outside of the boundary of the constraint).

\section*{Acknowledgments} \noindent The first author was supported by Ministry of Education and Science of Russian Federation in the frameworks of the project part of the state work quota (project number 1.3464.2017), the joint Taiwan NSC-Russia RFBR grant \textnumero 17-51-52022. 

\vskip0.2cm

\noindent We acknowledge useful discussions with Ivan Gudoshnikov (Mathematics, UT Dallas) and Supun Samarakoon (Mathematics, A\&M University) concerning the proof of Lemma~\ref{extra1}. The counter-example disproving the inequality 
$\|{\rm proj}(u,C)-{\rm proj}(u,D)\|\le d_H(C,D)$ (see Remark~\ref{remivan} and Fig.~\ref{contrfigure}) is also due to Ivan Gudoshnikov.




\end{document}